\documentclass[12pt]{amsart}

\usepackage{amsmath,amscd,amsfonts,amsthm,amssymb,eucal}
\usepackage{url}
\usepackage{fullpage}

\usepackage[]{latexsym}

\newtheorem{thm}{Theorem}[section]
\newtheorem{cor}[thm]{Corollary}
\newtheorem{lem}[thm]{Lemma}
\newtheorem{prop}[thm]{Proposition}

\newtheorem{conj}[thm]{Conjecture}

\theoremstyle{definition}
\newtheorem*{defn}{Definition}
\newtheorem*{rem}{Remark}
\newtheorem*{exmp}{Example}
\newtheorem*{ques}{Question}

\newenvironment{pf}{\par\noindent{\bf Proof.}\enspace\ignorespaces}{\qed\par\par}
\def\qed{\relax\ifmmode\hskip2em \Box\else\unskip\nobreak\hskip1em $\Box$\fi}

\newcommand{\red}{\mbox{red}}

\newcommand{\bQ}{{\mathbb{Q}}}

\newcommand{\bC}{{\mathbb{C}}}
\newcommand{\bP}{{\mathbb{P}}}
\newcommand{\bR}{{\mathbb{R}}}
\newcommand{\bH}{{\mathbb{H}}}
\newcommand{\bZ}{{\mathbb{Z}}}

\newcommand{\bG}{{\mathbb{G}}}

\newcommand{\cC}{{\mathbf{C}}}
\newcommand{\caC}{{\mathcal{C}}}
\newcommand{\cO}{{\mathcal{O}}}

\newcommand{\GL}{\mbox{GL}}
\newcommand{\SL}{\mbox{SL}}
\newcommand{\Sp}{\mbox{Sp}}
\newcommand{\Jac}{\mbox{Jac}}
\newcommand{\Aut}{\mbox{Aut}}
\newcommand{\End}{\mbox{End}}
\newcommand{\Pic}{\mbox{Pic}}
\newcommand{\Br}{\mbox{Br}}
\newcommand{\diam}{\mbox{diam}}

\newcommand{\plim}{\displaystyle{\lim_{\stackrel{\longleftarrow}{n}}}\,}

\begin{document}

\title{Trivial points on towers of curves}
\author{Xavier Xarles}
\address{Departament de Matem\`atiques\\Universitat Aut\`onoma de
Barcelona\\08193 Bellaterra, Barcelona, Catalonia}
\email{xarles@mat.uab.cat}

\thanks{The author was partially supported by the grant
MTM2009-10359.} \maketitle

\section{Introduction}

When studding the solutions of families of diophantine equations,
there are some (usually called) trivial solutions. Moreover, one
would like to know if for all members of the family with, may be, a
finite number of excepcions, the only solutions are the trivial
ones. In this paper we would like to study what should be a trivial
solution for the case that the family forms a tower. The idea is
that the trivial solutions are the solutions that are always there,
so they should be points that exists at all the levels of the tower.
And our main goal will be to find conditions in order to show when
there is a finite number of such trivial points, and also when there
are uniform bounds for the level of the tower where all the points
are trivial.

The main example to have in mind is the case of the modular towers:
consider the curves $X_1(p^n)$ for some fixed prime $p$, and $n\ge
0$. Then the trivial points should correspond to the cuspidal
points. And, moreover, there is a constant $N(p,d)$ depending on $p$
and $d\ge 1$ such that, for any number field $K$ with $[K:\bQ]\le
d$, and for any $n\ge N(p,d)$, the only $K$-rational points of
$X_1(p^n)$ are the trivial ones (see for example \cite{Ma}).

Another example was considered by the author in \cite{Xa}; the
result is similar to the case of modular curves. This paper is, in
some sense, a sequel of that paper where we investigate under
which circumstances one can get this type of results.

There are other recent papers treating similar problems. For
example in the paper \cite{CT} it is studied cases similar to the
modular towers, related to the inverse Galois problema (see also
\cite{Fri} for other generalizations), and in \cite{EHK} it is
studied also such (vary general) cases but for families indexed by
prime numbers.

The paper is organized as follows. In section 2 we introduce the
towers and their trivial points, and give some elementary results.
In section 3 and 4 we study the special cases of towers with genus
0 and 1, giving some partial results. In section 5 we give a
criterium for proving the finiteness of trivial points of a tower.
In section 6 we recall the well-known relation between the
unbounded gonality and the existence of uniform bounds (see
Theorem ??). The rest of the sections are dedicated to the
distinct methods to bound the gonality of a tower: geometric
methods, reduction modulo primes and counting points, and to
methods related to graphs.

This paper contains results concerning towers of curves that the
author collected during some years. The content was explained during
the ``Cuartas Jornadas de Teor\'{\i}a de N\'{u}meros'' in Bilbao, in July
2011. I would like to thank the organizers for the invitation to
give a talk, which motivated me to write this paper. I thank also
Enrique Gonz\'{a}lez Jim\'{e}nez, Joan Carles Lario, Francesc Bars, Pete
Clark, and Brian Conrad, for some conversations related to the
subject. I specially thank Bjorn Poonen for answering some doubts
concerning the gonality and for the help proving some results in
section 3.

\section{Notations, generalities and examples}

Given a field $K$, we will denote by $\bar K$ a fixed separable
closure of $K$ and by $G_K$ the absolute Galois group of $K$,
equal to $\Aut_K(\bar K)$, the automorphisms of the field $\bar K$
fixing the elements of $K$. Given any scheme $V$ over $K$, we will
denote by $\bar V$ the base change of $V$ to $\bar K$.

\begin{defn} Let $K$ be a field. By a tower of curves over $K$ we
mean a couple of sequences $\cC:=(\{C_n\}_{n\ge 0},
\{\varphi_n\}_{n\ge 1})$, where, for any $n\ge 0$, a positive
integer, $C_n$ are smooth projective algebraic curves defined over
$K$ and geometrically connected, and, for any strictly positive
integer $n>0$, the $\varphi_n:C_n \to C_{n-1}$ are non constant
morphisms as algebraic curves, of degree $>1$ for all $n\ge 0$.

If $m>n$, we will denote by $\varphi_{m,n}$ the morphism from $C_m$
to $C_n$ obtained composing the morphisms from $\varphi_m$ to
$\varphi_{n+1}$. We will denote also $\varphi_{n,n}$ the identity
map on $C_n$.
\end{defn}

\begin{defn} Given two towers of curves $\cC=(\{C_n\}_{n\ge
0},\{\varphi_n\}_{n\ge 1})$ and $\cC'=(\{C'_n\}_{n\ge
0},\{\varphi_n'\}_{n\ge 1})$, a morphism $\Psi$ from $\cC$ to
$\cC'$ is a collection of morphisms $\psi_n:C_n\to C'_{k_n}$ such
that $k_n\le k_m$ and
$\varphi'_{k_n,k_m}\psi_n=\psi_m\varphi_{n,m}$ if $n\le m$. A
special case is when $\cC$ is a subtower of $\cC'$, i.e
$C_n=C'_{k_n}$ and $\psi_n$ is the identity, for some progression
$k_0<k_1<k_2<\dots$.
\end{defn}

Given a smooth projective algebraic curve $C$ over a field $K$, we
will denote by $g(C)$ the genus of $C$, and by $\gamma(C)$ the
gonality of $C$ over the field $K$ (see section 5 for details).

\begin{defn}
Given a tower of curves $\cC$, we define the genus $g(\cC)$ to be
$\limsup g(C_n)$, and the gonality $\gamma(\cC)$ as $\limsup
\gamma(C_n)$.
\end{defn}

Given a tower of curves, we have by Hurwitz's theorem that
$g(C_n)\ge g(C_{n-1})$ for any $n$, so $0\le \limsup g(C_n)=\lim
g(C_n) \le +\infty$. The same fact for the gonality is not so
easy; one can find a proof of this result for example in
Proposition A1 in \cite{Po} (see also section 6).

\begin{lem} For any tower $\cC$, the genus can only be $g(\cC)=0$,
$1$ or $\infty$.
\end{lem}

\begin{pf} This is again due to Hurwtiz: if there exists a curve
$C_n$ in the tower with genus $g(C_n)>1$, then $g(C_{m+1})>g(C_m)$
for all $m\ge n$, so $g(\cC)=\infty$.
\end{pf}

\begin{exmp} One can construct genus 0 towers easily by
fixing rational functions $f_n(x)\in K(x)$ of degree $>1$. Over
algebraically closed fields, all genus 0 towers are of this type
(see section 3).

Also, to construct genus 1 towers, one can fix elliptic curves $E_n$
and isogenies $\varphi_n: E_n\to E_{n-1}$. An easy example is given
when $E_n=E$ and $\varphi_n$ is multiplication by some fixed integer
number for all $n\ge 0$ (see section 4).
\end{exmp}

\begin{defn} Given a tower of curves $\cC=(\{C_n\}_{n\ge 0},\{\varphi_n\}_{n\ge 1})$ defined over $K$,
the $K$-trivial points of $\cC$ in the level $n\ge 0$ are
$$ \cC(K)_n:=\{P\in C_n(K) \; | \; \forall m\ge n \ \exists P_m\in C_m(K)\,
\mbox{ such that } \varphi_{m,n}(P_m)=P\}.$$

If $d\ge 1$ is an integer, the  $K$-trivial points of $\cC$ in the
level $n\ge 0$ and degree $d$ are
$$ \cC^{(d)}(K)_n:=\bigcup_{L\subset \bar K \, , [L:K]\le d} \cC(L)_n .$$

Finally, the trivial points of $\cC$ in level $n$ are
$$ \cC^{(\infty)}(K)_n:= \bigcup_{d\ge 1}\cC^{(d)}(K)_n.$$
So the trivial points in level $n$ are the points that are
$L$-trivial in level $n$ for some finite extension $L/K$.

In the case of level 0 we will frequently omit it from the
notation.
\end{defn}

Observe that over an algebraically closed field the trivial points
of a tower of curves in the level $n$ is equal to all the rational
points of the curve $C_n$. So these definitions are only
interesting in fields such as number fields, or, more generally,
fields finitely generated over the prime fields.

Observe also that the $K$-rational trivial points in level $n$ are
equal to the imatge in $C_n(K)$ of the natural map $$\plim
C_n(K)\to C_n(K),$$ where the projective limit is taken with
respect to the maps $\varphi_n$.

The main problem we are going to analyze is when there are a finite
number of $K$-trivial points, and when a finite number of trivial
points. The next lemma is easy and we leave the proof to the reader.

\begin{lem} Let $\cC=(\{C_n\}_{n\ge 0},\{\varphi_n\}_{n\ge 1})$ be a tower of curves
over some field $K$.
\begin{enumerate}
\item For any $d=0,1,\dots,\infty$, the natural map $\varphi_{n,m}:\cC^{(d)}(K)_n\to \cC^{(d)}(K)_m$ is
surjective.
\item For any $d=0,1,\dots,\infty$, $\cC^{(d)}(K)_n \subseteq \varphi_{n,m}^{-1}(\cC^{(d)}(K)_m)$.
\item Two isomorphic towers have natural bijections between their
trivial points.
\item If $\cC'$ is a subtower of $\cC$, i.e. $C'_n=C_{k_n}$ and $\psi_n$ is the identity, for some
progression $k_0<k_1<k_2<\dots$, then
$\cC'^{(d)}(K)_n=\cC^{(d)}(K)_{k_n}$ for all $n\ge 0$ and all
$d=0,1,\dots,\infty$.
\end{enumerate}
\end{lem}

We end this section by given some specific examples of towers of
curves.

\begin{exmp} Consider a prime number $p$. We will denote by the $p$ Fermat
tower the tower with curves $C_n$ given by planar homogeneous
equations $X_0^{p^n}+X_1^{p^n}=X_2^{p^n}$, and maps
$\varphi_n(X_0,X_1,X_2)=[X_0^p:X_1^p:X_2^p]$. One can show that
the only trivial points of the tower are the ``trivial solutions"
$[a_0:a_1:a_2]$ with $a_0a_1a_2=0$.
\end{exmp}

\begin{exmp} Consider a prime number $p$. We will denote by the $p$
modular tower the tower with curves $C_n:=X_1(p^n)$ and natural
maps $\varphi_n$. One can show (see section 5) that the only
trivial points of the tower are the cuspidal points.
\end{exmp}

\begin{exmp} Consider the homogeneous polynomial $f_0:=X_0^2+X_1^2-X_2^2$, and
$C_n\subset \bP^{n+2}$ be the curves defined over $\bQ$ as the
zero set of the polynomials $f_0,f_1,\dots,f_n$, where
$f_n:=f_0(X_{n},X_{n+1},X_{n+2})$. For any field $K$, the
$K$-rational points of $C_n$ are in bijection with the Fibonacci
type sequences of squares of length $n+2$, that is sequences
$\{a_0,a_1,a_2,\dots,a_n\}$ such that $a_{n+2}=a_{n+1}+a_{n}$ and
all the elements $a_i$ are squares in $K$. We will call this tower
the Fibonacci tower.

Observe that we have four points $[\pm1,\, 0,\,\pm1,\, 1]\in
C_1(\bQ)$. One can show that $C_1$ is isomorphic to the elliptic
curve $E$ with Cremona Reference 32a2, and that $E(\bQ)$ has only
four points. So $C_1(\bQ)$ is form by that four points, thus
$C_2(\bQ)=\emptyset$ and, hence, $\cC(K)_0=\emptyset$.

Using the results in section 6 one can show that the degree $2$
points over $\bQ$ of the curve $C_2$, which has genus 5 and
gonality 4 (see section 7), inject inside the jacobian
$\Jac(C_2)$. Using results as in \cite{GX,GX2}, one can show that
the jacobian is isogenous to the product of 5 elliptic curves,
with Cremona references 32a2,32a2,48a1,96a1 and 96b1. All of them
have rank 0 and four rational points, so $\Jac(C_2)$ is finite and
computable. Using this result one gets that
$$C_2(\bQ)^{(2)}=\{[\sqrt{-1},\, \pm1,\, 0,\,\pm1,\, \pm1], [
\pm1,\, 0,\,\pm1,\, \pm1, \sqrt{2}]\},$$ which can be use to show
that $C_3(\bQ)^{(2)}=\emptyset$ and hence
$\cC(K)_0^{(2)}=\emptyset$. We do not know if
$\cC(K)_0^{(d)}=\emptyset$ for some $d\ge3$ (but we do know these
sets are finite using the results in section 6) or if
$\cC(K)_0^{\infty}=\emptyset$, or even finite.
\end{exmp}

\begin{conj} The Fibonacci tower has no trivial points over
any number field. Hence, the curves given by the system of
equations
$$X_0^2+X_1^2=X_2^2\  , \  X_1^2+X_2^2=X_3^2 ,\  \dots \   ,
X_n^2+X_{n+1}^2=X_{n+2}^2$$ inside $\bP^{n+2}$ have no rational
points for any number field $K$ and for $n$ large enough (in terms
of the degree of $K/\bQ$).
\end{conj}

\section{Genus 0 Towers}

Consider a tower of curves $\cC$ with genus 0. It is well known that
a genus 0 curve is either isomorphic to the projective line $\bP^1$
(and if and only if it has a rational point in your field), or
isomorphic to a conic curve (see for example Theorem A.4.3.1. in
\cite{HS}). In this second case, there exists some degree 2
extension of the field where the curve gets isomorphic to $\bP^1$
(and, even, there are infinitely many such degree 2 extensions if
the field is a number field).

\begin{defn} Given an enumerated set $\mathcal F:=\{f_n(x)\in K(x)\}_{n\ge 1}$ of rational
functions with degree $>1$, consider the tower of curves
$\cC_{\mathcal F}$ defined as $C_n:=\bP^1$ and $\varphi_n=f_n(x)$
for all $n\ge 0$. The special case that $f_n=f$ for all $n$ will
be denoted by $\cC_{f}$.\end{defn}

The next lemma was communicated to me by Bjorn Poonen.

\begin{lem} Let $C$ and $C'$ be genus 0 curves over a number
field $K$, and $f$ be a non-constant morphisms from  $C$ to $C'$.
Suppose that $C'(K)=\emptyset$. Then $C$ and $C'$ are isomorphic
(and the degree of $f$ is odd).
\end{lem}

\begin{pf} It is well known that the genus 0 curves $C$ without
$K$-rational points correspond to conics without points, so to
quaternion algebras over $K$, hence they give elements $x_C$ of
order 2 in the Brauer group $\Br(K)$ of $K$. The existence of the
map $f$ says us that $C$ has no rational points. We will see that
$\deg(f)x_C=x_{C'}$ in $\Br(K)$, hence $\deg(f)$ is odd and
$x_C=x_{C'}$, so $C$ is isomorphic to $C'$. Observe that we have a
natural map
$$ \bZ\cong \Pic(\bar C)^{G_K}\to H^1(G_K,\frac{\bar K(C)^*}{\bar K^*})\to H^2(G_K,\bar K^*)=\Br(K)$$
sending $1$ to $x_C$, given by the natural connecting homomorphisms,
which is functorial. On the other hand, the natural map $\Pic(\bar
C')^{G_K}\to \Pic(\bar C)^{G_K}$ is the multiplication by the degree
of $f$. Hence the result is deduced from the commutativity of the
natural diagram, which is easy.
\end{pf}

As a consequence of this lemma and the results cited above, we get
the following classification.

\begin{lem} Let $K$ be a number field, and let $\cC$ be a tower of
curves. Then
\begin{enumerate}
\item If $\cC(K)_0\ne \emptyset$, then there exists a set $\mathcal
F:=\{f_n(x)\in K(x)\}_{n\ge 1}$ of rational functions and an
isomorphism $\cC\cong \cC_{\mathcal F}$.
\item There always exists infinite degree 2 extensions $L/K$,
an enumerated set $\mathcal F:=\{f_n(x)\in L(x)\}_{n\ge 1}$ and an
isomorphism $\cC_L\cong \cC_{\mathcal F}$ defined over $L$.
\end{enumerate}
\end{lem}

\begin{pf} If $\cC(K)_0\ne \emptyset$, then there are points $P_n\in C_n(K)$ for all $n$,
hence all $C_n$ are isomorphic to $\bP^1$. But then the maps
$\varphi_n$ give us endomorphisms of $\bP^1$, so rational funcions
$f_n$. The same is true if $C_n(K)\ne \emptyset $ for all $n\ge
0$.

Now, suppose there exists $n\ge 0$ such that $C_n(K)=\emptyset$.
Using the previous lemma, we get that $C_m$ is isomorphic to $C_n$
for all $m\ge n$, and all isomorphic to a fixed conic $C$. Now we
only need to recall that for any fixed conic $C$ over a number
field there are an infinite number of quadratic extensions $L/K$
such that $C(L)\ne \emptyset$. \end{pf}

Now we are going to study the finiteness of $K$-trivial points.
Hence we can and will assume that $\cC(K)_0\ne \emptyset$. We will
only get some partial results, using the theory of heights.

\begin{rem} For the towers of the form $\cC_f$, with $f(x)\in K(x)$,
observe that $\cC(K)_n=\cC(K)_0$ for all $n$, and it contains the
set of periodic points of $f$: the points $P\in \bP^1(K)$ such that
$f^N(P)=P$ for some $N\ge 1$ (see \cite{HS}, B.4, for example).
\end{rem}

\begin{thm}\label{periodic} Let $K$ be a number field and let $f(x)\in K(x)$ be a rational
function of degree $d>1$. Then the set of $K$-rational trivial
points of level 0 of the tower $\cC_f$ is finite, and equal to the
set of periodic points of $f$.
\end{thm}
\begin{pf} We will show that the set of $K$-rational points is
equal to the set of periodic points, and finiteness will follow.
Consider the canonical height function $h_f$ associated to $f$
(see for example \cite{HS}, Theorem B.4.1.). Then
$h_f(f(P))=dh_f(P)$ for any $P\in \bP^1$, where $d>1$ is the
degree of $f$. If $P\in\cC_f(K)_0$, there exists $P_n\in \bP^1$
such that $f^n(P_n)=P$, so $h_f(P)=h_f(f_n(P_n))=d^nh_f(P_n)$.
Now, suppose that $h_f(P)\ne 0$. Then $h_f(P_n)=d^{-n}h_f(P)\le
h_f(P)$ are all distinct, so the set $\{P_n\ | \ n\ge 0\}$ is an
infinite set of points of $\bP ^1$ with bounded canonical height
$h_f$, hence with bounded absolute (logarithmic) height, which is
not possible (see \cite{HS}, theorem B.2.3). So $h_f(P)=0$ and
hence $P$ is preperiodic by Proposition B.4.2.(a) in \cite{HS}.
Now, if $P$ is preperiodic but not periodic, the set $\{P_n\ | \
n\ge 0\}$ is an infinite set of preperiodic points, which is again
no possible by Northcott theorem (Proposition B.4.2.(b) in
\cite{HS}).
\end{pf}

Observe that the finiteness of the set of $K$-rational trivial
points of degree $e>1$ will follow from the same result on the
periodic points, which is a conjecture (see for example conjecture
3.15 in \cite{SiD}).

\begin{ques} Are there genus 0 towers (of the form $\cC_{\mathcal
F}$) over a number field $K$ having an infinite number of
$K$-rational trivial points?
\end{ques}

Concerning the trivial points of the genus 0 towers, it is easy to
construct examples such that there are infinitely many of them,
and towers with only a finite number of them, as shown in the next
two examples.

\begin{exmp} Consider the special case $f(x)=x^2$ and $K=\bQ$. Then the set of
trivial points $\cC_{x^2}^{(\infty)}(K)_0$ of $\cC$ in level $0$
contains all the $n$-roots of unity for odd $\ge 1$ (and $x=0$ and
$\infty$), since they are defined in a finite extensi\'{o}n of $\bQ$ and
they are periodic for $f$. It is easy to show that they are all:
$$\cC_{x^2}^{(\infty)}(\bQ)_0=\{\xi\in\overline{\bQ}\ |\ \exists N\ge 1\mbox{ such that } \xi^N=1\}\cup\{0,\infty\}.$$
\end{exmp}

\begin{exmp} Take the genus 0 tower $\cC_{\mathcal F}$ defined by
$\mathcal F:=\{f_n(x)=x^{n+1}\in K(x)\}$, where $K$ is any number
field. Then the set of trivial points of $\cC$ in level $0$ is
equal to
$$\cC_{\mathcal F}^{(\infty)}(K)_0=\{0,1,\infty\}.$$
To show these, observe that $\alpha\in K$ is a $K$-rational
trivial point if and only if $\alpha$ has a $n$-th root for all
$n\ge 1$. But the only such numbers are $0$ and $1$ in any number
field. This last result can be shown using proving first that the
absolute logarithmic height of $\alpha$ must be $0$ (if $\alpha\ne
0$) as in the proof of the theorem \ref{periodic}, so $\alpha$
must be a root of unity. But the only root of unity which is and
$n$-th root of unity for all $n\ge1$ is $1$.
\end{exmp}

\section{Genus 1 Towers}

First of all, observe that, if we have a tower $\cC$ of genus 1
curves over a field $K$ such that there is a trivial point $P\in
\cC(K)_0$, we can use this point in order to get an explicit
description of the tower.

\begin{lem} Consider a tower $\cC=(\{C_n\}_{n\ge 0},\{\varphi_n\}_{n\ge 1})$ of genus 1
curves over a field $K$, and suppose there is a point $P\in
\cC(K)_0$.  Then the tower $\cC$ is isomorphic to a tower
$\mathcal E:= (\{E_n\}_{n\ge 0},\{\phi_n\}_{n\ge 1})$ where the
$E_n$ are elliptic curves and the $\phi_n: E_n\to E_{n-1}$ are
isogenies.
\end{lem}

\begin{pf} Let us fix a point $P_n\in C_n(K)$ such that
$\varphi_n(P_n)=P_{n-1}$. Consider the elliptic curve
$E_n:=\Jac(C_n)$, and the Abel-Jacobi map $\iota_n:C_n\to E_n$
given by the point $P_n\in C_n(K)$, which is an isomorphism of
curves. The maps $\phi_n:=\iota_{n-1}\varphi_n\iota_n^{-1}$ are
non-constant morphisms of curves between $E_n$ and $E_{n-1}$ which
send the identity point to the identity point. Hence they are
isogenies.
\end{pf}

\begin{defn} Given an elliptic curve $E$ over a field, denote by $\cO:=\End_K(E)$
the ring of endomorfisms of $E$ over $K$ as elliptic curve (if $E$
has no complex multiplication, and $K$ has characteristic 0, then
$\cO\cong \bZ$). Given a progression $\mathcal A:=\{a_n\in \cO | \
n=0,1,\dots \}$ of elements with degree $>1$, consider the tower
of curves $\cC_{E,\mathcal A}$ defined as $C_n:=E$ for all $n\ge
0$ and $\varphi_n$ equal to the $a_n$ for all $n\ge 1$. The
special case that $a_n=a$ for all $n$ will be denoted by
$\cC_{E,a}$.
\end{defn}

\begin{cor} Let $\cC=(\{C_n\}_{n\ge 0},\{\varphi_n\}_{n\ge 1})$ be a tower of genus 1
curves over a number field $K$ such that $\cC(K)_0\ne \emptyset$.
Then there is an elliptic curve $E$ defined over $K$ and a
progression $\mathcal A:=\{a_n\}$ such that the tower
$\cC_{E,\mathcal A}$ is isomorphic to a subtower of $\cC$.
\end{cor}

\begin{pf}
By applying the lemma we are reduced to the case that the $C_n$
are elliptic curves and the $\varphi_n$ are isogenies. A
well-known result of Faltings' (see \cite{Fal}) implies that there
is a finite number of elliptic curves isogenous to a given one.
Hence in the set of elliptic curves $C_n$, there are infinitely
many of them isomorphic to a given elliptic curve $E$. The result
is now easily deduced.\end{pf}

Now we can proof the finiteness of the $K$-rational trivial
points.

\begin{cor} Let $\cC=(\{C_n\}_{n\ge 0},\{\varphi_n\}_{n\ge 1})$ be a tower of genus 1
curves over a number field $K$. Then, for all $n\ge 0$ and all
$d\ge 1$, the set of $K$-trivial points of $\cC$ in the level
$n\ge 0$ and degree $d$ is finite.
\end{cor}
\begin{pf} Observe that it is sufficient to know the result for an isomorphic
subtower. Hence, by using the previous results, we are reduced to
the case that  $\cC=\cC_{E,\mathcal A}$ for some elliptic curve $E$
and some progression $\mathcal A:=\{a_n\}$ of natural numbers, or of
elements in a quadratic imaginary order in the CM case, and also the
case $n=0$.

Consider a point $P\in \cC^{(d)}(K)_0$. Let $L/K$ an extension of
degree $d$ such that $P\in \cC(L)_0$. Such an element in
$\cC(L)_0$ is a point $P\in E(L)$ that is divisible by
$b_n:=a_0a_1\cdots a_n$ for all $n$, hence by its norm to $\bZ$.
But $E(L)$ is finitely generated, so $P$ must be torsion. Thus
$\cC(L)_0\subset E(L)_{tors}$, the torsion subgroup, which is
finite, which proves the case $d=1$. In general, we get that
$$\cC^{(d)}(K)_0\subset \bigcup_{L\subset \bar K \, , [L:K]\le d}
E(L)_{tors}$$ which is again finite (a fact that can be deduced
from Merel's result \cite{Me}, or by using the results of section
6 due to Frey \cite{Fre} and the result of Abramovich bounding
below the gonality of the modular curves $X_1(N)$ \cite{Abr}, or,
even, the results in section 8).
\end{pf}

On the other hand, it is not true in general that the set of all
trivial points of a genus 1 tower is finite, as shown in this
example.

\begin{exmp}  The set of all
trivial points of the tower $\cC_{E,a}$, where $E$ is an elliptic
curve defined over a number field $K$ and $a>1$, is equal to the set
of torsion points of $E(\overline{K})$ with order prime with $a$:
$$\cC_{E,N}^{(\infty)}(K)=\{ P\in E(\bar K) \; | \; \exists m,
[m](P)=0 \mbox{ and } (m,a)=1 \}.$$
\end{exmp}

Finally, let us mention that there are genus 1 towers without
trivial points at all, as shown in the next example.

\begin{exmp} Let $E$ be an elliptic curve over a number field $K$.
Suppose that the Galois cohomology group $H^1(K,E)$ contains a
divisible element $\psi$, or, even less, an element divisible for
all powers of a prime $p$. Now, consider elements $\psi_n \in
H^1(K,E)$ such that $p\psi_n=\psi_{n-1}$ and $\psi_0=\psi$.

Recall that any element $\xi \in H^1(K,E)$ corresponds to a twist
$C_{\xi}$ of the curve $E$, that is, a genus 1 curve isomorphic to
$E$ over the algebraic clausure of $K$ (and, hence, with jacobian
isomorphic to $E$ over $K$). Moreover, the multiplication-by-$m$
in the group $H^1(K,E)$ correspons to a map $\phi_{m\xi}$ between
$C_{\xi}$ and $C_{m\xi}$ such that gives the multiplication-by-$m$
map between the corresponding jacobians (see for example \cite{LT}
for all this facts).

So we have a tower given by the curves $C_{\psi_n}$ and the maps
$\varphi_n:=\phi_{p\psi_n}$. Now, the elements $\psi_n$ have order
divisible by $p^n$, hence have index also divisible by $p^n$ (see
proposition 5 in \cite{LT}), which implies that the curves
$C_{\psi_n}$ do not have rational points in any extension with
degree $<p^n$. Hence the result.

Now, we only need to show the existence of such elliptic curves
$E$. But in fact, showing the existence of one such $E$ over $\bQ$
is sufficient. And this is known: take, for example, any elliptic
curve $E$ over $\bQ$ with finite number of $\bQ$-rational points
(see for example theorem D in \cite{CS}).
\end{exmp}

\section{Finiteness of trivial points and reduction}

Now we are going to consider towers with genus $\infty$, or,
equivalently, towers such that there is a curve $C_n$ with genus
$>1$. In this case, and when $K$ is a number field, the finitness
of the $K$-rational trivial points is clear, since, by Faltings'
theorem \cite{Fal}, the number of points in $C_n(K)$ is finite. So
we are mainly interested in the whole trivial points. Next example
will show that there are towers with genus infinite and an
infinite number of trivial points.

\begin{exmp} For $n\ge 0$, let $C_n$ be the smooth hyperelliptic curve defined
over $\bQ$ by the hyperelliptic equation $y^2=x^{2^n}-1$, and
consider the degree two maps $\varphi_n$ defined in the affine part
by $\varphi_n(x,y)=(x^2,y)$.

Now, take $\xi\in\overline{\bQ}$ a root of unity of odd degree, so
there exists an odd $N\ge 1$ such that $\xi^N=1$. Consider the field
$K_{\xi}$ generated by $\sqrt{\xi^i-1}$, for $i=1,\dots,N-1$; it is
a finite extension of $\bQ$ and, clearly, $(\xi,\sqrt{\xi-1})\in
C_0(K_{\xi})$ is a $K_{\xi}$-rational trivial point of the tower.

So, the trivial points of the tower include all the points of this
form, which are infinite.
\end{exmp}

Observe that hyperelliptic curves have an infinite number of
points of degree $2$, a result that generalizes to points of
degree $d$ and gonality $d$.

\begin{ques} For $d>1$, are there towers with infinite genus
over a number field $K$ having an infinite number of trivial points
of degree $d$?
\end{ques}

We will see in the next section that the answer of the question is
no when the gonality of the tower is infinite. But before we will
give a criterion for a tower to have a finite number of trivial
points.

\begin{defn} Let $\cC$ be a tower of curves over a number field
$K$ and a ring of integers $\cO$ (a Dedekind domain, not a field,
and with field of fractions $K$). By a proper model of $\cC$ over
$\cO$ we mean a collection of proper models $\caC_n$ of $C_n$ and
morphisms $\varphi_{n,\cO}:\caC_n\to \caC_{n-1}$ such that the
generic fiber is $\varphi_{n,\cO}\otimes_{\cO}K=\varphi_{n}$. We
will denote by $\widetilde{\varphi_{n,\wp}}:\caC_{n,\wp}\to
\caC_{n-1,\wp}$ the reduction of the morphism modulo some prime
$\wp$ of $\cO$ (and we will suppress the $\wp$ in the notation if
it is clear from the context).
\end{defn}

Observe that for any prime $\wp$ of $\cO$ we have a reduction map
$\red_{\wp}:\caC_n(\cO)\to \caC_{n,\wp}(k_{\wp})$, where $k_{\wp}$
is the residue field $\cO/\wp$. For any $P\in C_n(K)=\caC_n(\cO)$,
we have
$\red_{\wp}(\varphi_{n}(P))=\widetilde{\varphi_{n}}(\red_{\wp}(P)).$

\begin{thm}\label{gonalityred} Let $\cC$ be a tower of curves over a number field $K$ such that
$\Omega:=\cC(K)_n$ is finite for some $n\ge 0$. Fix a proper model
of $\cC$ over a ring of integers $\cO$ of $K$, and suppose that
for any prime $\wp$ of $\cO$ outside a finite number of primes,
there exists $m:=m_{\wp}\ge n$ such that
$C_m(k_{\wp})=\red_{\wp}(\varphi_{m,n}^{-1}(\Omega))$, where
$k_{\wp}$ is the residue field modulo $\wp$ and $\red_{\wp}$ is
the reduction map. Then $\Omega$ is the set of all trivial points
$\cC^{(\infty)}(K)_n$ of $\cC$, and there is a finite number.
\end{thm}

\begin{pf}
Let $L/K$ be a finite extension of $K$, and let $\cO_L$ be the
ring of integers of $L$. We are going to show that the set of
$L$-rational trivial points $\cC(L)_n$ is equal to $\Omega$.
Suppose in the contrary that there is a point
$P\in(\cC(L)_n\setminus \Omega)$.

Consider a prime ideal $\wp_L$ of $\cO_L$ such that there is a
prime ideal $\wp$ of $K$ divisible by $\wp_L$ and with equal
residue fields $k_{\wp_L}=k_{\wp}$ (there are an infinite number
of them). First, we show that then $\red_{\wp_L}(P)\in
\red_{\wp}(\Omega)$. Take $m:=m_{\wp}$ and $P_m\in C_m(L)$ such
that $\varphi_{m,n}(P_m)=P$. By hypothesis, $\red_{\wp_L}(P_m)\in
\red_{\wp}(\varphi_{m,n}^{-1}(\Omega))$. Hence
$$\widetilde{\varphi_{m,n}}(\red_{\wp_L}(P_m))=\red_{\wp_L}(\varphi_{m,n}(P_m))=\red_{\wp_L}(P)
\in
\widetilde{\varphi_{m,n}}(\red_{\wp}(\varphi_{m,n}^{-1}(\Omega)))=\red_{\wp}(\Omega).$$

So, we have an infinite number of primes $\wp_L$ of $\cO_L$ such
that $\red_{\wp_L}(P)\in \red_{\wp_L}(\Omega)$. Since the set
$\Omega$ is finite, there should exists a point $Q\in \Omega$ and
an infinite number of primes $\wp_L$ of $\cO_L$ such that
$\red_{\wp_L}(P)=\red_{\wp_L}(Q)$. But, given a proper curve
$\caC$ over a Dedekind domain $\cO$ and points $P\ne Q\in
\caC(\cO)$, the number of primes $\wp$ of $\cO$ such that
$\red_{\wp}(P)=\red_{\wp}(Q)$ is finite. Hence $P=Q$.
\end{pf}

\begin{cor} Suppose that the tower $\cC$ has a model with good reduction outside a
finite number of primes $S$ of $\cO$. For any prime $\wp$ of good
reduction of $\cC$, consider the tower of curves
$\cC_{\wp}=(\{\caC_{n,\wp}\},\{\widetilde{\varphi_{n,\wp}}\})$,
reduction modulo $\wp$ of the model.
 Suppose that there exists $n$ such that for any prime $\wp$
 outside a finite set containing $S$,
 $\red_{\wp}(\cC(K)_n)=\cC_{\wp}(k_{\wp})_n$. Then
 $\cC^{(\infty)}(K)_n=\cC(K)_n$.
\end{cor}

One can use this result to show that for some modular towers the
trivial points are the cusps. The next example is a generalization
of the curves $X_1(p^n)$ to higher dimensional abelian varieties.

\begin{exmp}  Let $K$ be a number field and $U/k$ a smooth geometrically connected
algebraic curve over $K$.  Let $\mathcal{A}\to U$ be an abelian
scheme of dimension $g\geq 1$, defined over $K$. Given a prime
number $p$ and a $n\ge 1$, consider the possibly disconnected
curve $\mathcal{A}[p^n]\to U$ over $U$, and let $U_n$ be a
geometrically connected component of it, such that the
multiplication-by-$p$ maps $\mathcal{A}[p^n]\to
\mathcal{A}[p^{n-1}]$ give maps $\phi_n:U_n\to U_{n-1}$ for all
$n\ge 2$. Let $C_n$ be the desingularization of some
projectivization of $U_n$, together with the natural maps
$\varphi_n:U_n\to U_{n-1}$. Since the $U_n$ are smooth over $U$,
the points $U_n(L)$ can be seen inside $C_n(L)$, for any extension
$L/K$. Then the trivial points $\cC^{(\infty)}(K)_n$ of the tower
$\cC=(\{C_n\}_{n\ge 0},\{\varphi_n\}_{n\ge 1})$ at level $n$ are
contained in $C_n(\bar K)\setminus U_n(\bar K)$, which is well
known to be finite.

This result can be shown by observing that the points of $s\in
U_n(L)$ correspond to a special fiber ${\mathcal A}_s$, which is
an abelian variety over $L$, together with a point $P\in {\mathcal
A}_s(L)$ of order exactly $p^n$. This point $s$ is a trivial point
if and only if for all $m\ge 1$, there exists a point $Q\in
{\mathcal A}_s(L)$ such that $[p^m]Q=P$; thus, $Q$ has exact order
$p^{n+m}$. Since the cardinal of the group of torsion points of an
abelian variety over a number field is finite (a fact that can be
proved by reducing modulo some primes $\ell$), there is no such a
$Q$ for $m\gg 1$. So there is no trivial point inside $U_n$.

Another way to show this result is by reducing modulo some prime
$\wp$ of $K$, such that the map $\mathcal{A}\to U$ has good
reduction, and does not divide $p$ (all primes except a finite
number of them verify these conditions). Then the tower $U_n$ has
good reduction at such a prime. The assertion is deduced from the
fact that the cardinal of the group of (torsion) points of an
abelian variety over a finite field is finite, which is trivial,
and then applying the corollary above.
\end{exmp}

\section{Genus $\infty$ Towers, trivial points and gonality}

Recall that the gonality $\gamma_K(C)$ of a curve $C$ over a field
$K$ is the minimum $m$ such that there exists a morphism
$\phi:C\to \bP^1$ of degree $m$ defined over $K$. For example,
hyperelliptic curves have gonality $2$. In the next proposition we
recall some properties of the gonality (see \cite{Po} for the
proofs).

\begin{prop} Let $K$ be any field, and let $C$ be an smooth and
projective curve with genus $g>1$ y gonality $\gamma_K(C)$. Then
\begin{enumerate}
\item $\gamma_K(C) \le 2g-2$.
\item If $C(K)\ne \emptyset$, then $\gamma_K(C)\le g$.
\item If $K=\bar K$ is algebraically closed, then $\gamma_K(C)\le \left\lfloor
\frac{g+3}{2}\right\rfloor $.
\item If $L/K$ is a field extension, then $\gamma_K(C)\ge\gamma_L(C)$.
\item If $K$ is a perfect field, $L/K$ is an algebraic field extension, $\gamma_L(C)>2$
and $C(K)\ne \emptyset$ then $\gamma_K(C)\le(\gamma_L(C)-1)^2$.
\item If $f:C\to C'$ is a non-constant $K$-morphism then
$\gamma_K(C)\le \deg(f)\gamma_K(C')$ and $\gamma_K(C')\le
\gamma_K(C)$.
\end{enumerate}
\end{prop}

The main tool we will use to relate the gonality with the
finiteness of trivial points is the following criterion of Frey
\cite{Fre}, proved also by Abramovich in his thesis, which is an
application of the main result of Faltings in \cite{Fal2}.

\begin{thm}[Frey] \label{Frey}
Let $C$ a curve over a number field $K$, with gonality $\gamma>2$
over $K$. Fix an algebraic closure $\overline{K}$ of $K$ and
consider the points of degree $d$ of $C$,
$$C^d(K):=\bigcup_{[L:K]\le d} C(L) \ \subset C(\overline{K})$$
where the union is over all the finite extensions of $K$ inside
$\overline{K}$ of degree $\le d$. Suppose that $2d<\gamma$. Then
$C^d(K)$ is finite.
\end{thm}

Using this criterion we can show the following result.

\begin{cor}\label{finite}
Let $\cC$ be a tower with infinite gonality. Then, for all $d\ge 1$
and $n\ge 0$, the set of $K$-trivial points of $\cC$ in the level
$n\ge 0$ and degree $d$ is finite.

Moreover, for all $d\ge 1$ there exists a constant $n_d$ (depending
on the tower $\cC$) such that, for any extension $L/K$ of degree
$\le d$, and for any $n\ge n_d$, $C_n(L)\subset \cC^{(d)}(K)_n$.
\end{cor}

\begin{cor}\label{finite2} If the set $\cC^{(\infty)}(K)_n$ is finite for some $n$,
and $\gamma(\cC)=\infty$, then for any extension $L/K$ of degree
$\le d$ and for any $n\ge n_d$, $C_n(L)\subset
\cC^{(\infty)}(K)_n$.
\end{cor}

\begin{ques} Are there towers with infinite gonality
over a number field $K$ having an infinite number of trivial points?
\end{ques}

\section{Gonality in the complete intersection case}

\begin{defn} For any $n\ge 1$, let $f_n(X_0,X_1,\dots,X_{n+1})$ be a
homogeneous polynomial of degree $d_n>1$. Consider the curves
$C_{n-1}\subset \bP^{n+1}$ defined as the zero set of the
polynomials $f_1,f_2,\dots, f_n$. We have a natural map $\varphi_n:
C_n\to C_{n-1}$ given by forgetting the last coordinate:
$\varphi_n([x_0:x_1:\dots:x_{n+2}])=[x_0:x_1:\dots:x_{n+1}]$. If the
polynomials $f_n$ are sufficiently general, then the curves $C_n$
are smooth and complete intersection, and we get a tower of curves.
We call these type of towers complete intersection towers.
\end{defn}

The main result of this section is that the complete intersection
towers have infinite gonality. The main tool is the following
theorem by Lazardsfeld (see Exercise 4.12. in \cite{La}), which is
a generalization of the well-known fact that a planar curve (so
given by an smooth projective model inside $\bP^2$) has gonality
equal to the degree of the model minus 1 if the curve has a
rational point.

\begin{thm}[Lazarsfeld] Let $C\subset \bP^n$ be a smooth complete
intersection of hypersurfaces of degrees $2\le a_1 \le a_2 \le
\cdots \le a_{r-1}$ over $\bC$. Then $\gamma(C) \ge (a_1-1)a_2\cdots
a_{r-1}$.
\end{thm}

\begin{cor} Any complete intersection
tower of curves has infinite gonality over any characteristic zero
field. \end{cor}

\begin{pf} Let $a_1:=\min{d_n}$ and let $n_0$ such that $d_{n_0}=a_1$.
Then the gonality of the curves $C_{n-1}$ for $n\ge n_0$ will be
bounded below by $$\gamma(C_{n-1})\ge \left(1-\frac
1{a_1}\right)\left(\prod_{i=1}^nd_n\right)$$ hence its limit goes to
infinite.
\end{pf}

Observe that, if the degree $d_1$ is the minimum of all the $d_n$,
and $C_0(K)\ne \emptyset$, then the gonality of $C_n$ over $K$ is
in fact equal to $(d_1-1)\left(\prod_{i=2}^nd_n\right)$, since the
map $\varphi_{n,1}:C_n\to C_0$ has degree
$\left(\prod_{i=2}^nd_n\right)$, which composed with the map of
degree $d_1-1$ from $C_0$ to $\bP^1$ given by the rational point
has the desired degree.

\begin{exmp} Fix an homogeneous irreducible polynomial $f(X_0,X_1,X_2)$ of degree
$d>1$ and defined over a number field. Suppose that the curve
projective $C_0 \subset \bP^{2}$ defined as the zero set of the
polynomials $f$ is non singular and geometrically connected.
Consider now the complete intersection tower $\cC_{f}$ of curves
$C_{n-1}\subset \bP^{n+1}$ defined as the zero set of the
polynomials $f_1:=f,f_2:=f(X_1,X_2,X_3),\dots,
f_n:=f(X_{n-2},X_{n-1},X_n)$. Then the curve $C_n$ has gonality
$\ge d^n-d^{n-1}$, so the tower has infinite gonality, with
equality exactly when $C_0(K)\ne \emptyset$.

Hence, if $K$ is a number field, and by Corollary \ref{finite}, we
get that the curves $C_n$ have only the trivial points over a
finite extension $L/K$ for $n$ large enough, depending only on the
tower and the degree $[L:K]$.

And, in case we know the set of trivial points is finite and
computable (for example, by theorem \ref{gonalityred}, we get that
the points of the curves $C_n$  over a finite extension $L/K$ for
$n$ large enough, depending only on $\cC$ and the degree $[L:K]$,
by Corollary \ref{finite2}.
\end{exmp}

One can also use the theorem for other type of towers of curves, a
generalization of the $p^n$ Fermat tower.

\begin{exmp} Fix an homogeneous irreducible polynomial $f_0(X_0,X_1,X_2)$ of degree
$d>1$ and defined over a number field, and a progression
$\{a_n\}_{n\ge 1}$ of integers $a_n\ge 2$. Suppose that the curve
projective $C_0 \subset \bP^{2}$ defined as the zero set of the
polynomials $f_0$ is non singular and geometrically connected.
Consider now the complete intersection tower $\cC$ of curves
$C_n\subset \bP^{2}$ defined as the zero set of the polynomial
$f_n:=f_{n-1}(X_0^{a_n},X_1^{a_n},X_2^{a_n})$, and the maps
$\varphi_n$ given by
$\varphi_n(X_0,X_1,X_2)=[X_0^{a_n}:X_1^{a_n}:X_2^{a_n}]$. Then the
curves $C_n$ are planar curves of degree $da_1a_2\cdots a_n$, and
hence with gonality $\gamma(C_n)\ge da_1a_2\cdots a_n-1\ge 2^n-1$.
\end{exmp}

\section{Gonality and reduction}

The following proposition was shown in \cite{Xa}, proposition 5,
and it will allow us to bound by below the gonality by just
counting points ``modulo p".

\begin{prop}\label{rgonality} Let $C$ be a curve over a number
field, and let $\wp$ be a prime of good reduction of the curve,
with residue field $k_{\wp}$. Denote by $C'$ the reduction of the
curve $C$ modulo $\wp$. Then the gonality $\gamma(C)$ of $C$
satisfies that
$$\gamma(C_K) \ge \gamma(C'_{k_{\wp}})\ge \frac{\sharp C'(k_{\wp})}{\sharp k_{\wp}+1}.$$
\end{prop}

Observe that the proposition has two parts: first, that the
gonality does not increase under reduction modulo a prime. Second,
that the gonality over a finite field is bounded below by the
number of points.

\begin{defn} Let $K$ be a number field and let $\cC$ be a tower of
curves over$K$. We will say that a prime $\wp$ of the ring of
integers $\cO$ of $K$ is a prime of good reduction of the tower if
there exists a proper model of $\cC$ over $\cO_{(\wp)}$, the
localization of $\cO$ at $\wp$, such that $\caC_{n,\wp}$ are
smooth and projective curves, and the maps
$\widetilde{\varphi_{n,\wp}}$ are non constant.
\end{defn}

\begin{cor} Let $\cC=(\{C_n\}_{n\ge 0},\{\varphi_n\}_{n\ge 1})$ be a tower of
curves over a number field $K$. Suppose that there exists a prime
$\wp$ of good reduction of the tower, and suppose that
$$\lim_{n\to \infty} \sharp  \caC_{n,\wp}(k_{\wp}) =+\infty.$$
Then the tower $\cC$ has infinite gonality.
\end{cor}

\begin{pf} Using the proposition \ref{rgonality} one gets that
$$\gamma(C_n) \ge \gamma(\caC_{n,\wp})\ge \frac{\sharp \caC_{n,\wp}(k_{\wp})}{\sharp k_{\wp}+1},$$
hence the result. \end{pf}

\section{Gonality and Cayley-Schreier graphs}

In this section we are going to follow the ideas originating in the
work of Zograf \cite{Zog} and Abramovich \cite{Abr}, and developed
in a recent paper of Ellenberg, Hall and Kowalski \cite{EHK}. The
idea is to show that certain \'{e}tale towers of (possibly affine)
curves have infinite gonality if the associated Cayley-Schreier
graphs form an expanding family (or, more generally, verify some
growing condition in the first non-trivial eigenvalue of the
combinatorial laplacian operator).

Suppose that we have a tower of curves $\cC$ defined over a number
field such that the maps $\varphi_{n,0}$ are \'{e}tale (i.e.
non-ramified) outside a fixed Zarisky closed set $Z$. So, we have
open subsets $U_i=C_i\setminus Z$ of $C_i$, together with maps
$\varphi_n:U_n\to U_{n-1}$, which are \'{e}tale, and the original tower
$\cC$ is obtained projectivizing the curves $U_i$ (the case
$U_i=C_i$ is also considered).

Fix, for all $i\ge 0$, a point $x_i$ in $U_i(\bar K)$ such that
$\varphi_i(x_i)=x_{i-1}$, and a generating set $S$ of the
topological fundamental group $G:=\pi_1({U_0}_{\bC},x_0)$. Consider
the Cayley-Schreier graphs $\Gamma_i=C(N_i,S)$ associated to the
finite quotient sets
$$N_i:=G/H_i=\pi_1({U_0}_{\bC},x_0)/\pi_1({U_i}_{\bC},x_i).$$

Recall that the graphs $\Gamma_i=C(G/H_i,S)$ have vertex set
$V(\Gamma_i)=G/H_i$, and with (possibly multiple) edges from vertex
$xH_i$ to vertex $sxH_i$ for all $s \in S$; hence, they are
$r$-regular graphs for $r=|S|$.

Define the combinatorial Laplacian operator of a $r$-regular graph
$\Gamma$ as $rId-A(\Gamma_i)$, where $A(\Gamma)$ is the adjacency
matrix of $\Gamma$. We compute the eigenvalues of $\Gamma$, which
are positive real numbers, and let $\lambda_1(G)$ to be the
smallest non-zero of them.

Observe that there are maps of graphs $\Gamma_i\to \Gamma_{i-1}$ for
all $i\ge 1$, and that such maps are unramified: the preimatge of
any vertex is formed always by $k$ vertices, where $k$, the degree
of the map, is fixed.

\begin{defn} A tower of graphs is a couple $(\{\Gamma_i\},\{\phi_i\})$
where $\Gamma_i$ are graphs for any $i\ge 0$ and $\phi_i:\Gamma_i\to
\Gamma_{i-1}$ are surjective maps of graphs. We say that the tower
is unramified if all the maps $\phi_i$ for $i\ge 1$ are unramified.
\end{defn}

Following ideas from the paper \cite{EHK}, we will bound the
gonality of the curves $C_i$ by imposing some condition on the first
non-zero eigenvalue $\lambda_1(\Gamma_i)$  of the combinatorial
laplacian operator on $\Gamma_i$.

\begin{thm} Let $\cC=(\{C_n\},\{\varphi_n\})$ be a tower of curves defined over $\bC$
such that the maps $\varphi_{n,0}$ are \'{e}tale outside a
fixed Zarisky closed set $Z$, and with the genus $g(C_0)>1$.
Consider the open subsets $U_i=C_i\setminus Z$ of $C_i$, together
with maps $\varphi_n:U_n\to U_{n-1}$ and points $x_i$ in $U_i(\bar
K)$ such that $\varphi_i(x_i)=x_{i-1}$. Let $\{\Gamma_i,\phi_i\}$ be
the unramified tower of Cayley-Schreier graphs
$C(\pi_1({U_0}_{\bC},x_0)/\pi_1({U_i}_{\bC},x_i),S)$. Suppose that
$\lim_{i\to \infty} \lambda_1(\Gamma_i)|V(\Gamma_i)|=\infty$. Then
$\lim_{i\to \infty}  \gamma(C_i)=\infty$.
\end{thm}

\begin{pf}
First of all, observe that, by the Hurwitz's formula
$$g(C_i)-1\ge \deg(\varphi_{i,0})(g(C_0)-1)\ge \deg(\varphi_{i,0})$$
(where we have equality exactly if $\phi_{i,0}$ are unramified).
Now, the degree of $\varphi_{i,0}$ is exactly equal to the index
of $\pi_1({U_i}_{\bC},x_i)$ inside $\pi_1({U_0}_{\bC},x_0)$, which
is equal to the number of vertices of $\Gamma_i$. So we get that
$g(C_i)-1\ge |V(\Gamma_i)|$.

Now, we will find a formula relating the $\lambda_1(\Gamma_i)$ to
the gonality of $C_i$ and $g(C_i)$. We will follow the proof of
Theorem 8 (b) in \cite{EHK}, and we will only sketch the proof.
Since the genus of $C_i$ is $>1$ for all $i$, we can write $U_i$ as
$G_i\backslash \bH$ for some discrete subgroup of $PSL_2(\bR)$. The
hyperbolic area $\mu_i(U_i)$ is then finite and the Poincar\'{e} metric
induces a Laplacian operator on the $L^2$-space. Following Li and
Yau \cite{LY}, and Abramovich \cite{Abr}, one has that
$$\gamma(C_i)\ge \frac 1{8\pi}\lambda_1(U_i)\mu(U_i)$$
where $\lambda_1(U_i)$ is the first non-trivial eigenvalue of the
laplacian operator $-div(grad)$.

Now, using the Gauss-Bonnet theorem, one gets
$$\mu(U_i)=-2\pi\chi(U_i)\ge -2\pi\chi(C_i)=-4\pi(1-g(C_i)).$$
Using the comparison principle of Brooks \cite{Bro} and Burger
\cite{Bur}, one gets that there exists a constant $c>0$, depending
only on $U_0$ and on $S$, such that
$$\lambda_1(U_i)\ge c \lambda_1(\Gamma_i).$$
Hence, combining all the results, that
$$\gamma(U_i)\ge 2c \lambda_1(\Gamma_i)(g(C_i)-1)\ge  2c \lambda_1(\Gamma_i)|V(\Gamma_i)|$$
and hence the result. \end{pf}

We say that a family of graphs is an expander if $\lim_{i\to
\infty}|\Gamma_i|=\infty$ and $\lambda_1(\Gamma_i)\ge c$ for some
constant $c$. We will say that it is esperantist if there exists
some constant $A\ge 0$ such that
$$\lambda_1(\Gamma_i)\ge \frac c{(\log (2|\Gamma_i|))^A}.$$

Observe that, if our family (in fact, tower) of Cayley-Schreier
graphs $\{\Gamma_i\}$ is an expander (or it is esperantist), then
they verify a fortiori the hypothesis of the theorem. Hence, the
following constructions give towers of curves $\cC$ defined over a
number field $K$ with infinite gonality. Consider $U_0$ a smooth
geometrically connected algebraic curve over a number field $K$, and
suppose we have an epimorfism of groups $p:\pi_1(U_0(\bC),x_0)\to
G$, where the group $G$ is one of the cases below. Take finite-index
subgroups $H_n$ of $G$ such that $H_n\varsubsetneq H_{n-1}$ and
$H_0=G$, and consider the \'{e}tale coverings $U_n\to U_0$ associated to
the subgroups $p^{-1}(H)$. Finally, consider the projectivizations
and desingularizations $C_n$ of these curves $U_n$.

\begin{enumerate}
\item If $G$ is a finite-index subgroup in ${\mathbf G}(\bQ)\cap \GL_m(\bZ)$, where
${\mathbf G} \subset {\mathbf \GL}_m$ is a semisimple algebraic
subgroup, defined over $\bQ$, and $G$ has real rank at least 2
(for example, $G$ can be a finite-index subgroup of $\SL_n(\bZ)$,
$n > 3$, or of $\Sp_{2g}(\bZ)$, $g > 2$) and $S$ is an arbitrary
finite set of generators of $G$, and $H_n$ arbitrary normal
subgroups (by property (T) of Kazhdan, see \cite{BHV}).
\item If $G$ is a subgroup of $\SL_n(\bZ)$ which is Zarisky-dense in
$\SL_d$, for $d>1$, $S$ is an arbitrary finite set of generators of
$G$, $p$ is a prime number sufficiently large (depending on $G$) and
$H_n=p^n \SL_d(\bZ)$, so $G/H_n\cong \SL_n(\bZ/p^n\bZ)$, by
\cite{BG1} and \cite{BG2}.
\item If $G$ is any Zarisky dense subgroup of ${\mathbf G}/\bZ_p$,
an arbitrary split semisimple algebraic group, and we consider the
tower of Cayley graphs of ${\mathbf G}(\bZ/p^n\bZ)$ with respect
to any symmetric set of generators, by the results of Dinai
\cite{Di1},\cite{Di2} concerning the diameter of this graphs: the
diameter is less than $c\log(|{\mathbf G}(\bZ/p^n\bZ)|)^d$, for
some constants $c$ and $d$.
\end{enumerate}

To show some of these cases, one needs to know that there is a
relation between the  first non-trivial eigenvalue of the
combinatorial laplacian operator and the diameter (longest
shortest path between any to pair of vertices) for any (regular)
graph $\Gamma$. For example, Diaconis and Saloff-Coste \cite{DSC}
showed that
$$\lambda_1(C(G,S))\ge \frac{1}{|S|\diam(C(G,S))^2},$$
if $C(G,S)$ is a Cayley Graph associated to a finite group $G$
with symmetric set of generators $S$. Hence, if e have normal
subgroups $H_i\unlhd G$ inside a group $G$, with $H_i\subseteq
H_{i-1}$, the hypothesis of the theorem is verified if
$$\lim_{i\to \infty} \frac{|G/H_i|}{\diam(C(G/H_i,S))^2}=+\infty$$

On the other hand, it is easy to construct towers of curves not
verifying the growing condition on the first non-trivial
eigenvalue of the combinatorial laplacian operator (and having
bounded gonality), as in the following trivial example.

\begin{exmp} Consider the tower of curves $\cC_{x^2}$, so with $C_n=\bP^1$
and maps given by $\varphi_n(x)=x^2$. Of course this tower has no
infinite gonality. These maps $\varphi_n$ are unramified outside
$x=0$ and $\infty$, so they give unramified selfmaps of
$U_i=\bG_m=\bP^1\setminus\{0,\infty\}$. Choose, for example,
$x_n:=1\in C_n(\bQ)$. Then one has that
$\pi_1({U_i}_{\bC},x_i)\cong \bZ$, and for the set $S:=\{1,-1\}$,
we get that the Cayley-Schreier graphs
$$C(\pi_1({U_0}_{\bC},x_0)/\pi_1({U_i}_{\bC},x_i),S)=C(\bZ/2^n\bZ,\{\pm1\})=\Gamma_{2^n}$$
where $\Gamma_n$ denotes the cycle graph form by a cycle with $n$
vertices and $n$ edges. It is well known that the eigenvalues of
the combinatorial laplacian operator for these cycle graphs are
$\lambda_k(\Gamma_n):=2-2\cos(2k\pi/n)$ for $k=0,\dots, n-1$.
Hence,
$$\lim_{i\to \infty} \lambda_1(\Gamma_{2^i})|V(\Gamma_{2^i})|=\lim_{i\to \infty} 2^i(2-2\cos(\frac
{\pi}{2^{i-1}}))=0.$$
\end{exmp}

Observe that there is another completely different relation of the
gonality with graph theory, developed by Baker and Norine
\cite{BN} and specially by Baker in \cite{Ba}.

Suppose $K$ is a field complete with respect to a discrete
valuation, and let $C$ be a curve over $K$ having a regular
semistable model over the ring of integers of $K$. Consider
$\Gamma$ the dual graph of the reduction of $X$, where the
vertexes are the irreducible components, and the edges correspond
to the intersection points. Then, there is a notion of gonality
for a finite graph and the gonality of $C$ is bounded below by the
gonality of $\Gamma$ (see Corollary 3.2 in \cite{Ba}). Using this
result it is not difficult to construct towers of (Mumford) curves
over $\bQ_p$ having infinite gonality.

All these results are also related to tropicalizations of
algebraic curves, and how to bound the gonality from the
tropicalization, a subject that we will explore in the future.

\end{document}